\documentclass[a4]{article}
\pdfoutput=1

\usepackage{pgf}
\usepackage{amsfonts}
\usepackage{fullpage}
\usepackage{tikz}
\usepackage{gnuplot-lua-tikz}

\begin{document}
  \title{On the effects of scaling on the performance of Ipopt}

\author{ J.~D.~Hogg and J.~A.~Scott\thanks{This work was supported by EPSRC grant EP/I013067/1}\\
Scientific Computing Department, \\
STFC Rutherford Appleton Laboratory,\\
Harwell Oxford, Didcot OX11 0QX, UK}
   
   \maketitle
   \begin{abstract}
      The open-source nonlinear solver Ipopt
      ({\tt https://projects.coin-or.org/Ipopt}) is a widely-used software package for
      the solution of large-scale non-linear optimization problems. At its heart, it employs a third-party
      linear solver to solve a series of sparse symmetric indefinite systems. The speed, 
      accuracy and robustness of the chosen linear solver is critical to the overall performance of
      Ipopt.       In some instances, it can be beneficial to scale the linear system
      before it is solved.

      In this paper, different scaling algorithms
      are employed within Ipopt with a new linear solver {\tt HSL\_MA97}
      from the HSL mathematical software library
      ({\tt http://www.hsl.rl.ac.uk/}). An extensive collection of
      problems from the CUTEr test set ({\tt http://www.cuter.rl.ac.uk/})
      is used to illustrate the effects of scaling.
   \end{abstract}

\section{Introduction}

   The Ipopt \cite{wabi:2006} package is designed to solve large-scale non-linear
   optimization problems of the form
   \begin{eqnarray*}
      \min_{x\in\mathbb{R}^n} && f(x) \\
      \mathrm{s.t.}  && g(x)=0, \\
      && x_L \le x \le x_U,
   \end{eqnarray*}
   where $f(x): \mathbb{R}^n \rightarrow \mathbb{R}$ is the objective function
   and $g(x): \mathbb{R}^n \rightarrow \mathbb{R}^m$ are the constraint
   functions. Note that more general formulations can be expressed in this form.
   Ipopt implements a primal-dual interior point filter line search
   algorithm for large-scale nonlinear programming. This involves
   solving the indefinite sparse symmetric linear system
   \begin{equation} \label{aug system}
       \left(\begin{array}{cc}
         W_k + \Sigma_k + \delta_w I & A_k \\
         A_k^T & -\delta_c I
      \end{array}\right)\left(\begin{array}{c}
         d_k^x \\
         d_k^\lambda
      \end{array}\right) = \left(\begin{array}{c}
         b_k^x \\
         b_k^\lambda
      \end{array}\right)
   \end{equation}
   where $\delta_w$ and $\delta_c$ are chosen such that the matrix has
   inertia $(n,m,0)$. $W_k$ and $A_k$ are the Hessian of the Lagrangian and the
   Jacobian of $g(x)$, respectively, evaluated at the current trial point.
   $\Sigma_k$ is a positive semi-definite diagonal matrix that has values approaching $0$ and
   $+\infty$ as the algorithm converges to a solution. The immediate result
   of this is that the condition number of the system may increase dramatically
   as optimality is approached. Further, it is critical that the solution method for
   (\ref{aug system}) accurately reports the inertia of the matrix to allow
   the Ipopt algorithm to choose $\delta_w$ and $\delta_c$.
   The Ipopt algorithm also uses multiple second-order corrections computed
   using the same system matrix as in (\ref{aug system}) with different
   right-hand sides, so efficient reuse of information from the initial solution
   is desirable.

   For these reasons, a sparse direct method is used to solve the system.
   If, for convenience of notation, we denote the indefinite system
   (\ref{aug system}) by 
   \begin{equation} \label{Ax=b}
      Ay = b,
   \end{equation}
    a direct method computes a factorization of the form
   $PAP^T=LDL^T$, where $P$ is a permutation matrix
    (or, more generally, a product of permutation matrices) that holds the pivot order
and is chosen to try to preserve sparsity and limit growth in the size of the
    factor entries,
    $L$ is a unit lower triangular matrix and $D$ is a block diagonal 
    matrix with blocks of order $1$ or $2$. The solution process is completed by 
   performing forward and back substitutions (that is, by first solving a lower 
   triangular system, a system involving $D$ and then an upper triangular
   system). The quality of the solution is
   generally dependent on a stability threshold pivoting parameter $u$: the larger $u$ is, the more accurate the solution
   but this can be at the cost of additional floating-point operations and more fill in $L$
   (see, for example, Duff, Erisman and Reid \cite{duer:86}, page 98).
   For some factorizations (particularly those where a small $u$ was used), it
   may be necessary to use a refinement process to improve the quality of the computed
   solution; normally this is iterative refinement, though more sophisticated
   schemes are sometimes used \cite{Arioli:2008,hosc:2010a}.

   In some applications, matrix scaling can be used to reduce both  the computational time and 
   memory requirements of the direct solver while also improving the accuracy of the computed solution,
   thus limiting the need for iterative refinement and also reducing the number of iterations used by Ipopt.
   We can apply a symmetric (diagonal) scaling matrix $S$ to (\ref{Ax=b})
   \begin{displaymath}
      SASS^{-1}y = Sb.
   \end{displaymath}
   This is equivalent to solving the following system
   \begin{displaymath}
      \hat{A}z  =  \hat{b},  
       y  =  Sz,
   \end{displaymath}
where  $\hat{A} = SAS$  and $ \hat{b} = Sb$.
The Ipopt ``Hints and Tips'' webpage states ``{\it
If you have trouble finding a solution with Ipopt, it might sometimes 
help to increase the accuracy of the computation of the search directions,
 which is done by solving a linear system. There are a number of ways to do this:
 First, you can tell Ipopt to scale the linear system before it is sent to the linear solver}''. 
How to find a good scaling $S$ is still an open question, but a 
number of scalings for sparse linear systems have been proposed and are widely used in
 many different application areas. In this paper, we experiment 
with scalings available from the HSL mathematical software
library \cite{hsl:2011}. In particular, 
we examine the use of a number of different scaling algorithms 
with the new sparse solver {\tt HSL\_MA97}   \cite{hosc:2011a} within Ipopt. We use a large 
collection of problems from the CUTEr test set.
Our aim is to illustrate
the effectiveness of scaling but also to show how it can significantly add to the Ipopt
runtime. Our key contribution is a detailed study of the effects
of different scaling strategies used with Ipopt that will assist users
in making an informed decision about how and when scaling should be tried.
We also introduce new heuristics that can be used within Ipopt to limit how
often the scaling matrix $S$ is recomputed without increasing the
number of Ipopt iterations, thus reducing not only the scaling time but also
the total Ipopt runtime.

The rest of this paper is organised as follows. We end this
section by outlining the experimental setup we will use.
In Section~\ref{scaling}, we briefly describe the different scaling algorithms
that are implemented by routines in the HSL library and will be used in our study.
Numerical results for these scalings used with {\tt HSL\_MA97} within Ipopt
are presented in Section~\ref{effect scalings}. In Section~\ref{dynamic scaling},
we look at the effect of selectively applying scaling as the Ipopt
algorithm progresses. Our recommendations for users are summarised in
Section~\ref{advice} and some final remarks are made in
Section~\ref{conclusions}.

\subsection{Test environment} \label{expt setup}

   \begin{figure}
      \caption{ \label{mitchell}
         Description of test machine.      }
      \begin{minipage}{0.5\linewidth}
         \centering
%         \pgfimage[height=7cm]{fig/mitchell}
         \pgfimage[height=7cm]{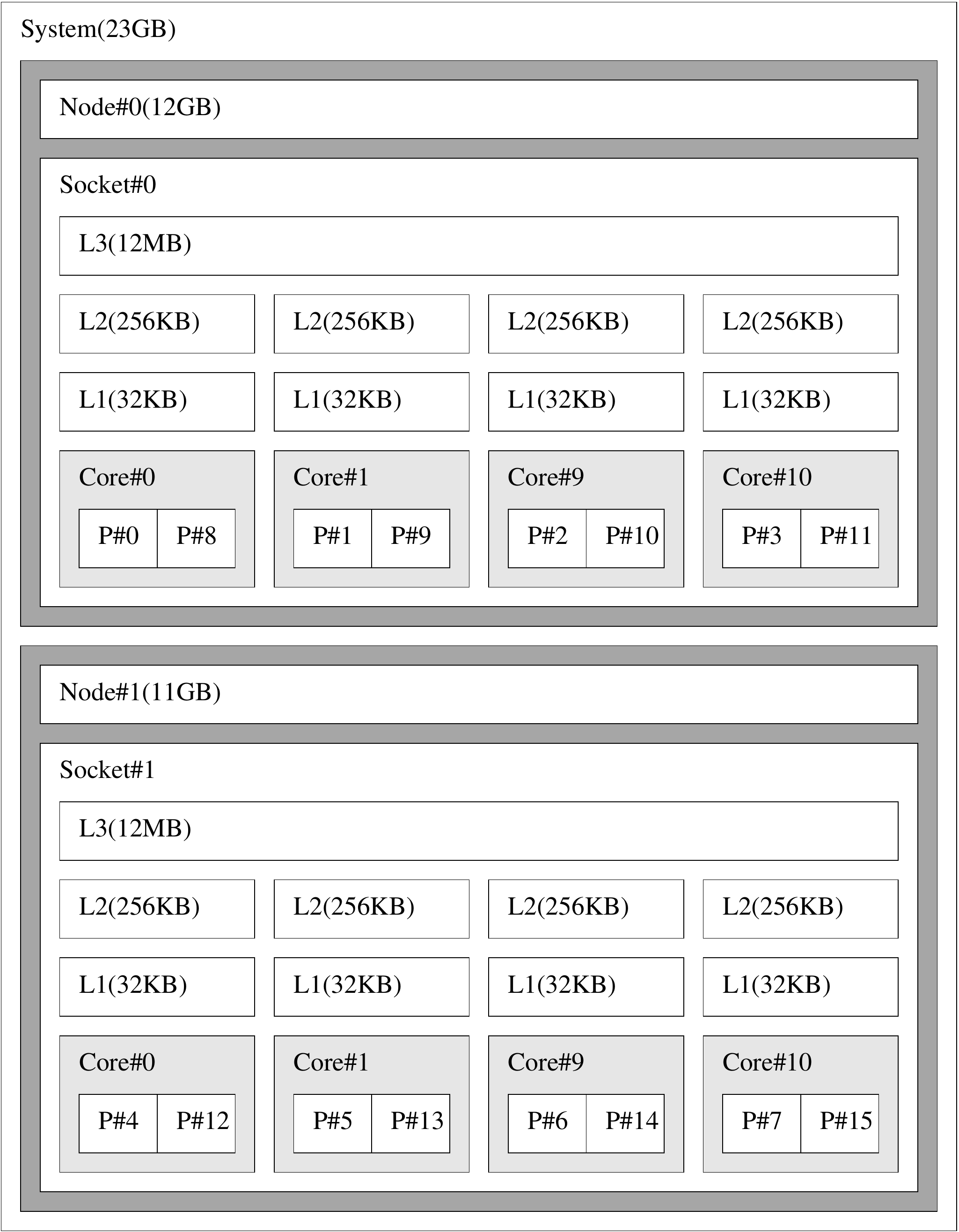}
      \end{minipage}
      \begin{minipage}{0.5\linewidth}
         \begin{tabular}{ll}
            \multicolumn{2}{c}{\textbf{Machine specification}} \\
            \hline
            \textbf{Processor} & 2 $\times$ Intel Xeon E5620\\
            \textbf{Physical Cores} & 8 \\
            \textbf{Memory} & 24 GB \\
            \textbf{Compiler} & Intel Fortran 12.0.0 \\
            & ifort -g -fast -openmp \\
            \textbf{BLAS} & MKL 10.3.0 \\
         \end{tabular}
      \end{minipage}
   \end{figure}

All our experiments are performed using {\tt HSL\_MA97} (Version 2.0.0) on the test machine summarised
   in Figure~\ref{mitchell}. Although {\tt HSL\_MA97} is a parallel code,
   all reported timings are serial times. We use Ipopt version 3.10. 
 %  The driver 
 %  for {\tt HSL\_MA97} is based on the Ipopt driver for the HSL solver
 %  {\tt HSL\_MA86}, but is adapted to support the full range of features offered by
 %  {\tt HSL\_MA97}.

   The CUTEr experiments are run in double precision using the Ipopt CUTEr
   interface and a development version of the CUTEr tools that supports
   64-bit architectures. Of the 1093 CUTEr problems, we discard 534 that take
   less than 0.1 seconds to solve using the HSL solver {\tt MA57} with no scaling ({\tt MA57}
   was chosen to select the test problems since it is our most widely used modern sparse
   direct solver). We then
   discard 173 problems that Ipopt fails to solve with any combination of
   (HSL) linear solver and scaling tested (note that although, in this paper,
   we only report on using {\tt HSL\_MA97}, we have performed extensive experiments
   with other HSL linear solvers and the results of these experiments 
  led to the set of discarded problems that we failed to solve). This leaves us with a set of
   386 non-trivial tractable problems, and it is to this set we refer to with
   the phrase ``all problems'' henceforth.
   Except as otherwise stated, the default Ipopt settings are used for {\tt HSL\_MA97} and
   CUTEr. In particular, an automatic choice between approximate minimum degree (AMD) ordering and 
   nested dissection (MeTiS) ordering is used. In addition, for each problem, the initial
   threshold pivoting parameter $u$ is set to $10^{-8}$. Note that this Ipopt choice is significantly smaller
   than the default setting within {\tt HSL\_MA97}  of $10^{-2}$; the use of small $u$ in an optimization
   context has been discussed, for example, in \cite{fome:93,giss:96}.
   For each test problem, a virtual memory limit of 24 GB and time limit of 1000 seconds is placed
   on each run. Any run exceeding these limits is classified as having failed.

%\subsection{Performance profiles for comparing scalings}
%
%   Because we have a large test set, many of our results are presented using performance profiles, 
%as described in \cite{Dolan:2002}.
%   Let $i \in \cal S$ be a scaling algorithm we wish to compare, with statistic (e.g. time) $s_{ij} \geq 0$ 
%for problem $j\in\cal T$, and that the smaller this statistic the better the scaling is considered to be. 
%If a scaling leads to a problem not being solved,   its statistic is considered to be $\infty$.
%
%   For $j\in \cal T$, let $\hat{s}_j = \min \{s_{ij}; i \in \cal S\}$ be the statistic for the 
%best performing scaling on that problem. Then for $\alpha \geq 1$ and each $i \in \cal S$
% we define the indicator function
%   \[
%      k(s_{ij},\hat{s}_{j},\alpha) = \left\{ \begin{array}{rl}
%         1 & \mbox{if} \;\; s_{ij} \leq \alpha \hat{s}_{j} \\ 0 & \mbox{otherwise.}
%      \end{array}\right.
%   \]
%   The performance profile of scaling $i$ is then given by the function
%   \[
%      p_i(\alpha) = \frac{\sum_{j \in \cal T} k(s_{ij},\hat{s}_{j},\alpha)}{|\cal T|}, \;\; \alpha \geq 1.
%   \]
%This can then be plotted on a graph showing all other scalings. The closer the line gets to the top left 
%corner, the better the scaling is. 

\section{The scalings}
\label{scaling}

\vspace{3mm}
\noindent
{\bf No Scaling}
\vspace{2mm}

\noindent
   Many problems can be solved without the use of scaling and, as our numerical experiments will show, this 
can often result in the fastest solution times (particularly for small problems). The default option within Ipopt 
(controlled by the parameter {\tt linear\_system\_scaling}) is not
to prescale the linear systems before calling the solver. However, it should be noted
that some of the linear solvers that are available for use within
Ipopt incorporate   scaling algorithms internally and these may be used by default.
   
\vspace{3mm}
\noindent
{\bf \texttt{MC19} and \texttt{MC30}}
\vspace{2mm}

\noindent  
 {\tt MC19} is the oldest scaling routine available in HSL; it is an implementation 
of the method described in the 1972 paper of Curtis and Reid~\cite{Curtis:1972}. 
Its aim is to make the entries of the scaled matrix close to one by minimizing the 
sum of the squares of the logarithms of the moduli of the entries.
That is, it scales $a_{ij}$ to
   \begin{displaymath}
      a_{ij} \exp(r_i + c_j),
   \end{displaymath}
   with the aim of minimizing the sum of the squares of the logarithms of the entries:
   \begin{eqnarray*}
       \min_{r,c} & \sum_{a_{ij}\ne0} \log |a_{ij}| + r_i + c_j .
   \end{eqnarray*}
   This is achieved by a specialized conjugate gradient algorithm. 
%In the symmetric case, the row
%and column scaling factors  are equal.

{\tt MC19} is available as part of the HSL Archive ({\tt http://www.hsl.rl.ac.uk/archive/}). 
This is a collection of older routines, many of
which have been superseded by routines in the main HSL library. In particular, {\tt MC19}
has been superseded by {\tt MC29} for unsymmetric matrices and by {\tt MC30}
for symmetric matrices. These codes are again based on \cite{Curtis:1972}.
The main differences between {\tt MC19} and {\tt MC30} are that the former uses
only single precision internally to compute the scaling factors, while the
latter uses double, and  {\tt MC30} uses a symmetric iteration, while
{\tt MC19} uses an unsymmetric iteration with the result averaged between rows
and columns at the end to obtain symmetric scaling factors.

\vspace{3mm}
\noindent
{\bf\texttt{MC64}}
\vspace{2mm}

\noindent 
Given an $n \times n$ matrix $A= \{a_{ij}\}$, the HSL package {\tt MC64} seeks to
find a matching of the rows to the columns such
that the product of the matched entries is maximised.
That is, it finds a permutation $\sigma= \{\sigma(i)\}$ 
of the integers 1 to $n$ that solves the maximum product matching problem
$$
   \max_{\sigma} \prod_{i=1}^n|a_{i\sigma(i)}|.
$$
The matrix entries
$a_{i\sigma(i)}$ corresponding to the solution $\sigma$ are said to be {\it matched}.
Once such a matching has been found, row and column scaling factors $r_i$ and
$c_{\sigma(i)}$ may be found such that all entries of the scaled matrix $S_rAS_c$, where
$S_r = \mathrm{diag}(r_i)$ and $S_c = \mathrm{diag}(c_{\sigma(i)})$, are less
than or equal to one in absolute value, with the matched entries having
absolute value one.

In the symmetric case, symmetry needs to be preserved.
Duff and Pralet~\cite{dupr:2005} start by treating $A$
as unsymmetric and compute the row and column scaling $r_i$ and $c_{\sigma(i)}$.
They show that if $A$ is non-singular, the geometric average of the
row and column scaling factors is sufficient to maintain desirable properties.
That is, 
they build a diagonal scaling matrix $S$ with entries
\[ 
   s_i = \sqrt{r_ic_i}
\]
 and compute the symmetrically scaled matrix $SAS$.
Again, the entries in the scaled matrix that are in the matching have absolute value one
while the rest have absolute value less than or equal to one.
The symmetrized scaling is included in the package {\tt HSL\_MC64},
a Fortran 95 version of the earlier code {\tt MC64}, with a number of
additional options. 

\vspace{3mm}
\noindent
{\bf {\tt MC77}}
\vspace{2mm}

\noindent
Equilibration is a particular form of scaling in which the rows and columns of the matrix are scaled
so that they have approximately the same norm. The HSL package
{\tt MC77} uses an iterative procedure \cite{Ruiz:2001} to attempt to make all row and 
column norms of the matrix unity for a user-specified geometric norm $\|\cdot\|_p$.  
The infinity norm is the default within {\tt MC77}. It produces a matrix whose rows 
and columns have maximum entry of exactly one and has good convergence properties. 
The one norm produces a matrix whose row and column sums are exactly one (a doubly 
stochastic matrix).

Recently, Ruiz and U\c{c}ar~\cite{ruuc:2011}  proposed combining the use of the infinity and one-norms.
Specifically, they  perform one step of  infinity norm scaling followed by
a few steps of one norm scaling (with more steps of the  one norm scaling used for their hard test problems).
This combination is the {\tt MC77} option available within {\tt HSL\_MA97},
with one step of  infinity norm scaling followed by
three of  one norm scaling.

\vspace{3mm}
\noindent
{\bf\texttt{HSL\_MC80}}
\vspace{2mm}

\noindent
If the matched entries returned by {\tt MC64} are unsymmetrically permuted onto
the diagonal, they provide good pivot candidates for an unsymmetric
factorization. To obtain similar benefits in a symmetric factorization, we
must instead symmetrically permute these entries on to the subdiagonals and
use $2\times2$ pivots.

Entry $a_{ij}$ is symmetrically permuted onto the subdiagonal by any
permutation in which $i$ and $j$ become adjacent in the final order.
It is not necessary to apply this condition for all matched entries. In fact,
it is sufficient if each row or column participates in a single $2\times2$
pivot. We therefore construct a subset of the matching such that an index
$k$ is only included in a pair once as either a row or column.
If this is done prior to the fill-reducing ordering step of the linear solver,
an alternative ordering can be found by considering the columns $i$ and $j$
in a pair as a single entity and applying standard reordering techniques. 
The resultant matching-based ordering is fill-reducing and has large
entries on the diagonal or subdiagonal for most rows and columns.

This approach is already implemented as part of the sparse solver
PARDISO~\cite{schenk:2004,scga:2006} and has been discussed in the context of
Ipopt by Schenk, W\"achter and Hagemann~\cite{scwh:2007}. The code
\texttt{HSL\_MC80} represents our implementation of their technique.

The analyse phase of the linear solver {\tt HSL\_MA97} optionally computes the
matching-based ordering and corresponding scaling by calling {\tt HSL\_MC80}.
Note that this approach has a number of potential downsides:
\begin{itemize}
   \item The cost of the analyse phase can be substantially increased.
   \item The analyse phase may need to be rerun whenever the numerical values
      change.
   \item The predicted fill-in will, in most cases, be higher than if numerical
      values were not taken into account when ordering.
\end{itemize}
However, for problems that require substantial modifications to the pivot
sequence during the numerical factorization to maintain stability, the actual
cost of the factorization may be reduced \cite{hosc:2012b}.
Note that the scaling used is the same as returned by {\tt MC64}, however the
ordering of the matrix is different.

\section{Effects of scaling on every iteration} \label{effect scalings}
In this section, we look at recomputing and applying the scaling at each Ipopt iteration.
Figure \ref{full scaling cmp} shows a performance profile \cite{Dolan:2002}
comparison of running Ipopt  with each of the  scalings discussed in 
Section~\ref{scaling} and with no scaling.
This demonstrates that not scaling leads to the fastest
solution times for most of the test problems. However, the test set of 386 problems
is skewed towards
smaller problems for which, as later results will demonstrate,
 the scaling time can represent a significant overhead
(and can dominate the linear solver time).
Figure \ref{large scaling cmp} shows the same profile restricted to the 54
problems that take more than 0.1s (in serial) per iteration. For these problems
(which we will refer to as the {\it large} problems) there is a much
weaker differentiation between scaling and not scaling. If we consider 
the asymptotes of these performance profiles, we obtain
a figure on reliability (the percentage of problems satisfactorily solved);
these are summarised in Table~\ref{scaling reliability}.
\begin{figure}
   \caption{ \label{full scaling cmp}
      Performance profile of times using different scalings (all problems).
   }
    \vspace{3mm}

   \centering
   \input{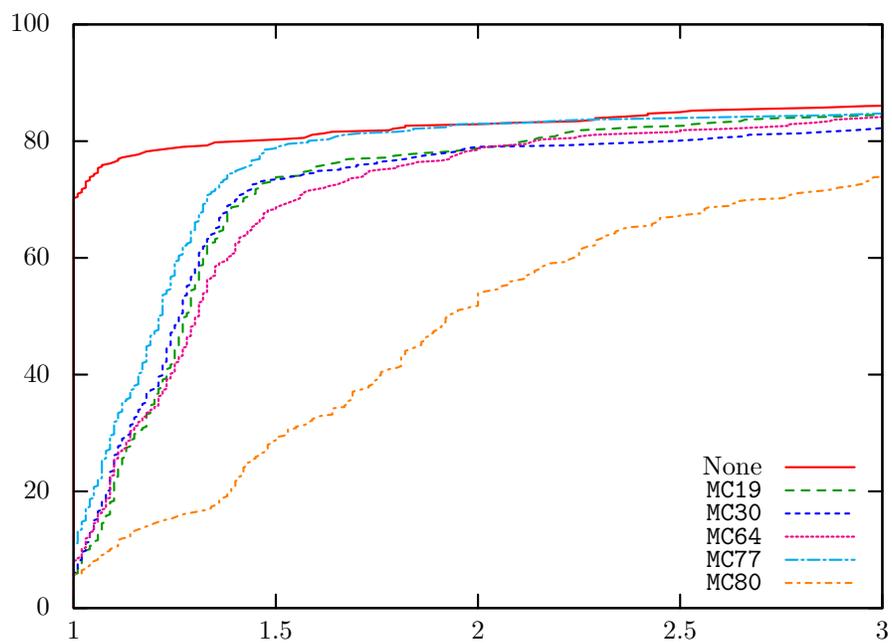}
\end{figure}
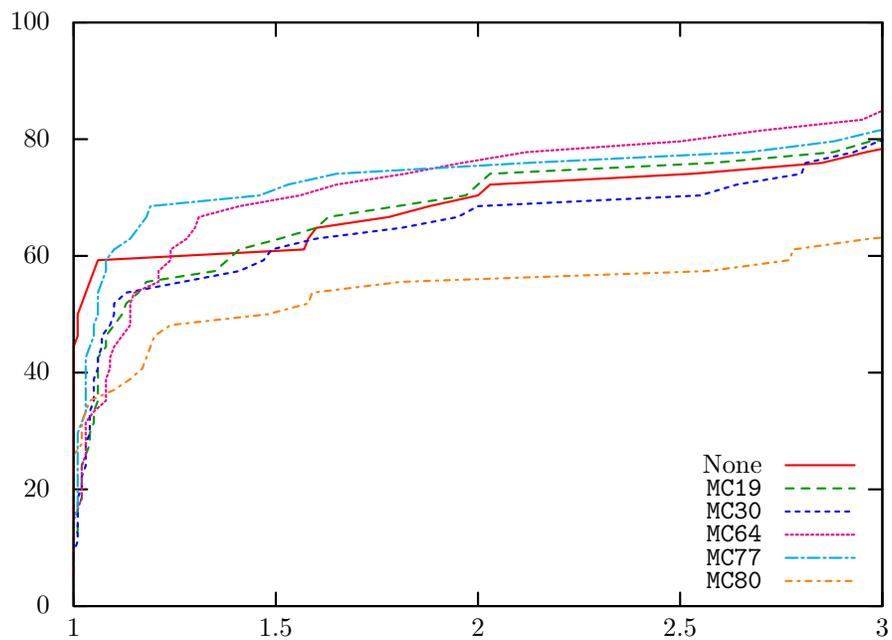
\begin{figure}
   \caption{ \label{large scaling cmp}
      Performance profile of times using different scalings (large problems).
   }
    \vspace{3mm}

   \centering
   \begin{tikzpicture}[gnuplot]
%% generated with GNUPLOT 4.4p3 (Lua 5.1.5; terminal rev. 97, script rev. 96a)
%% Mon 10 Dec 2012 10:45:36 GMT
\gpcolor{gp lt color border}
\gpsetlinetype{gp lt border}
\gpsetlinewidth{2.00}
\draw[gp path] (1.196,0.616)--(1.376,0.616);
\draw[gp path] (11.947,0.616)--(11.767,0.616);
\node[gp node right] at (1.012,0.616) { 0};
\draw[gp path] (1.196,2.169)--(1.376,2.169);
\draw[gp path] (11.947,2.169)--(11.767,2.169);
\node[gp node right] at (1.012,2.169) { 20};
\draw[gp path] (1.196,3.722)--(1.376,3.722);
\draw[gp path] (11.947,3.722)--(11.767,3.722);
\node[gp node right] at (1.012,3.722) { 40};
\draw[gp path] (1.196,5.275)--(1.376,5.275);
\draw[gp path] (11.947,5.275)--(11.767,5.275);
\node[gp node right] at (1.012,5.275) { 60};
\draw[gp path] (1.196,6.828)--(1.376,6.828);
\draw[gp path] (11.947,6.828)--(11.767,6.828);
\node[gp node right] at (1.012,6.828) { 80};
\draw[gp path] (1.196,8.381)--(1.376,8.381);
\draw[gp path] (11.947,8.381)--(11.767,8.381);
\node[gp node right] at (1.012,8.381) { 100};
\draw[gp path] (1.196,0.616)--(1.196,0.796);
\draw[gp path] (1.196,8.381)--(1.196,8.201);
\node[gp node center] at (1.196,0.308) { 1};
\draw[gp path] (3.884,0.616)--(3.884,0.796);
\draw[gp path] (3.884,8.381)--(3.884,8.201);
\node[gp node center] at (3.884,0.308) { 1.5};
\draw[gp path] (6.572,0.616)--(6.572,0.796);
\draw[gp path] (6.572,8.381)--(6.572,8.201);
\node[gp node center] at (6.572,0.308) { 2};
\draw[gp path] (9.259,0.616)--(9.259,0.796);
\draw[gp path] (9.259,8.381)--(9.259,8.201);
\node[gp node center] at (9.259,0.308) { 2.5};
\draw[gp path] (11.947,0.616)--(11.947,0.796);
\draw[gp path] (11.947,8.381)--(11.947,8.201);
\node[gp node center] at (11.947,0.308) { 3};
\draw[gp path] (1.196,8.381)--(1.196,0.616)--(11.947,0.616)--(11.947,8.381)--cycle;
\node[gp node right] at (10.479,2.490) {None};
\gpcolor{gp lt color 0}
\gpsetlinetype{gp lt plot 0}
\draw[gp path] (10.663,2.490)--(11.579,2.490);
\draw[gp path] (1.196,0.760)--(1.196,0.903)--(1.196,1.048)--(1.196,1.191)--(1.196,1.335)%
  --(1.196,1.479)--(1.196,1.622)--(1.196,1.766)--(1.196,1.910)--(1.196,2.054)--(1.196,2.198)%
  --(1.196,2.341)--(1.196,2.485)--(1.196,2.629)--(1.196,2.773)--(1.196,2.917)--(1.196,3.060)%
  --(1.196,3.204)--(1.196,3.349)--(1.196,3.492)--(1.196,3.636)--(1.196,3.779)--(1.196,3.923)%
  --(1.196,4.067)--(1.250,4.211)--(1.250,4.355)--(1.250,4.499)--(1.304,4.642)--(1.357,4.786)%
  --(1.411,4.930)--(1.465,5.074)--(1.519,5.218)--(4.260,5.361)--(4.314,5.505)--(4.421,5.648)%
  --(5.389,5.793)--(5.926,5.937)--(6.572,6.080)--(6.733,6.224)--(9.421,6.368)--(11.141,6.512)%
  --(11.732,6.656)--(11.947,6.700);
\gpcolor{gp lt color border}
\node[gp node right] at (10.479,2.182) {\tt MC19};
\gpcolor{gp lt color 1}
\gpsetlinetype{gp lt plot 1}
\draw[gp path] (10.663,2.182)--(11.579,2.182);
\draw[gp path] (1.196,0.760)--(1.196,0.903)--(1.196,1.048)--(1.196,1.191)--(1.196,1.335)%
  --(1.196,1.479)--(1.250,1.622)--(1.250,1.766)--(1.250,1.910)--(1.304,2.054)--(1.304,2.198)%
  --(1.304,2.341)--(1.304,2.485)--(1.357,2.629)--(1.411,2.773)--(1.411,2.917)--(1.465,3.060)%
  --(1.465,3.204)--(1.519,3.349)--(1.519,3.492)--(1.519,3.636)--(1.519,3.779)--(1.519,3.923)%
  --(1.626,4.067)--(1.626,4.211)--(1.734,4.355)--(1.841,4.499)--(1.895,4.642)--(2.056,4.786)%
  --(2.164,4.930)--(3.077,5.074)--(3.239,5.218)--(3.400,5.361)--(3.938,5.505)--(4.421,5.648)%
  --(4.583,5.793)--(5.496,5.937)--(6.410,6.080)--(6.572,6.224)--(6.733,6.368)--(9.689,6.512)%
  --(11.302,6.656)--(11.786,6.799)--(11.947,6.814);
\gpcolor{gp lt color border}
\node[gp node right] at (10.479,1.874) {\tt MC30};
\gpcolor{gp lt color 2}
\gpsetlinetype{gp lt plot 2}
\draw[gp path] (10.663,1.874)--(11.579,1.874);
\draw[gp path] (1.196,0.760)--(1.196,0.903)--(1.196,1.048)--(1.196,1.191)--(1.196,1.335)%
  --(1.250,1.479)--(1.250,1.622)--(1.250,1.766)--(1.250,1.910)--(1.250,2.054)--(1.304,2.198)%
  --(1.304,2.341)--(1.357,2.485)--(1.357,2.629)--(1.357,2.773)--(1.411,2.917)--(1.411,3.060)%
  --(1.411,3.204)--(1.465,3.349)--(1.465,3.492)--(1.465,3.636)--(1.519,3.779)--(1.519,3.923)%
  --(1.572,4.067)--(1.572,4.211)--(1.680,4.355)--(1.734,4.499)--(1.734,4.642)--(1.895,4.786)%
  --(2.647,4.930)--(3.400,5.074)--(3.722,5.218)--(3.830,5.361)--(4.421,5.505)--(5.550,5.648)%
  --(6.303,5.793)--(6.572,5.937)--(9.528,6.080)--(10.012,6.224)--(10.872,6.368)--(10.926,6.512)%
  --(11.571,6.656)--(11.893,6.799)--(11.947,6.808);
\gpcolor{gp lt color border}
\node[gp node right] at (10.479,1.566) {\tt MC64};
\gpcolor{gp lt color 3}
\gpsetlinetype{gp lt plot 3}
\draw[gp path] (10.663,1.566)--(11.579,1.566);
\draw[gp path] (1.196,0.760)--(1.196,0.903)--(1.196,1.048)--(1.196,1.191)--(1.196,1.335)%
  --(1.196,1.479)--(1.196,1.622)--(1.196,1.766)--(1.250,1.910)--(1.304,2.054)--(1.304,2.198)%
  --(1.304,2.341)--(1.304,2.485)--(1.357,2.629)--(1.357,2.773)--(1.357,2.917)--(1.357,3.060)%
  --(1.465,3.204)--(1.626,3.349)--(1.626,3.492)--(1.626,3.636)--(1.680,3.779)--(1.680,3.923)%
  --(1.734,4.067)--(1.841,4.211)--(1.949,4.355)--(1.949,4.499)--(1.949,4.642)--(2.002,4.786)%
  --(2.325,4.930)--(2.325,5.074)--(2.486,5.218)--(2.486,5.361)--(2.701,5.505)--(2.809,5.648)%
  --(2.862,5.793)--(3.400,5.937)--(4.206,6.080)--(4.690,6.224)--(5.604,6.368)--(6.356,6.512)%
  --(7.217,6.656)--(9.259,6.799)--(10.334,6.943)--(11.678,7.087)--(11.947,7.207);
\gpcolor{gp lt color border}
\node[gp node right] at (10.479,1.258) {\tt MC77};
\gpcolor{gp lt color 4}
\gpsetlinetype{gp lt plot 4}
\draw[gp path] (10.663,1.258)--(11.579,1.258);
\draw[gp path] (1.196,0.760)--(1.196,0.903)--(1.196,1.048)--(1.196,1.191)--(1.196,1.335)%
  --(1.196,1.479)--(1.196,1.622)--(1.196,1.766)--(1.250,1.910)--(1.250,2.054)--(1.250,2.198)%
  --(1.250,2.341)--(1.250,2.485)--(1.250,2.629)--(1.250,2.773)--(1.250,2.917)--(1.304,3.060)%
  --(1.357,3.204)--(1.357,3.349)--(1.357,3.492)--(1.357,3.636)--(1.357,3.779)--(1.357,3.923)%
  --(1.411,4.067)--(1.465,4.211)--(1.465,4.355)--(1.519,4.499)--(1.519,4.642)--(1.519,4.786)%
  --(1.572,4.930)--(1.626,5.074)--(1.626,5.218)--(1.734,5.361)--(1.949,5.505)--(2.056,5.648)%
  --(2.164,5.793)--(2.217,5.937)--(3.669,6.080)--(4.045,6.224)--(4.690,6.368)--(7.217,6.512)%
  --(10.173,6.656)--(11.302,6.799)--(11.893,6.943)--(11.947,6.948);
\gpcolor{gp lt color border}
\node[gp node right] at (10.479,0.950) {\tt MC80};
\gpcolor{gp lt color 5}
\gpsetlinetype{gp lt plot 5}
\draw[gp path] (10.663,0.950)--(11.579,0.950);
\draw[gp path] (1.196,0.760)--(1.196,0.903)--(1.196,1.048)--(1.196,1.191)--(1.196,1.335)%
  --(1.196,1.479)--(1.196,1.622)--(1.196,1.766)--(1.196,1.910)--(1.196,2.054)--(1.196,2.198)%
  --(1.196,2.341)--(1.196,2.485)--(1.196,2.629)--(1.304,2.773)--(1.304,2.917)--(1.304,3.060)%
  --(1.357,3.204)--(1.411,3.349)--(1.734,3.492)--(1.949,3.636)--(2.110,3.779)--(2.164,3.923)%
  --(2.217,4.067)--(2.271,4.211)--(2.486,4.355)--(3.776,4.499)--(4.314,4.642)--(4.368,4.786)%
  --(5.550,4.930)--(9.636,5.074)--(10.711,5.218)--(10.764,5.361)--(11.786,5.505)--(11.947,5.520);
\gpcolor{gp lt color border}
\gpsetlinetype{gp lt border}
\draw[gp path] (1.196,8.381)--(1.196,0.616)--(11.947,0.616)--(11.947,8.381)--cycle;
%% coordinates of the plot area
\gpdefrectangularnode{gp plot 1}{\pgfpoint{1.196cm}{0.616cm}}{\pgfpoint{11.947cm}{8.381cm}}
\end{tikzpicture}
\end{figure}
\begin{table}[htbp]
   \caption{ \label{scaling reliability}
      Reliability of different scalings used on the set of all (non-trivial, tractable) problems and the subset of large problems.
   }
    \vspace{3mm}

   \centering
   \small
   \begin{tabular}{l|cc}
      \hline
                     & All    & Large \\
                     & problems & problems \\
      \hline
      None           & 92.8\% & 94.4\% \\
      \texttt{MC19}  & 92.0\% & 90.7\% \\
      \texttt{MC30}  & 92.0\% & 92.6\% \\
      \texttt{MC64}  & 92.5\% & 96.3\% \\
      \texttt{MC77}  & 87.8\% & 90.7\% \\
      \texttt{MC80}  & 92.2\% & 90.7\% \\
      \hline
   \end{tabular}
\end{table}

\begin{figure}
   \caption{ \label{ndelay cmp}
      Performance profile of the average number of delayed pivots per factorization for different 
scalings (all problems).
   }
    \vspace{3mm}

   \centering
   \input{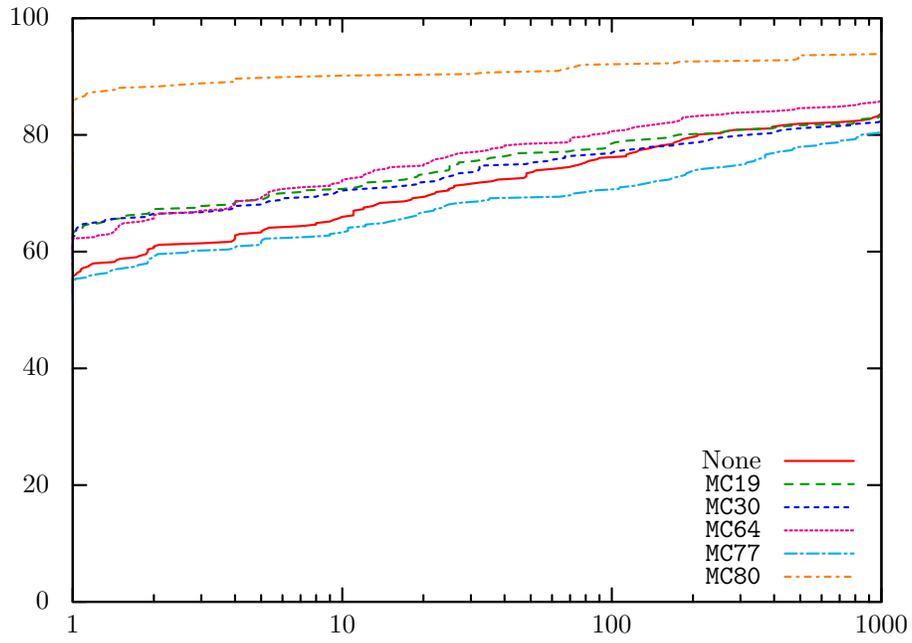}
\end{figure}

\begin{figure}
   \caption{ \label{nflops cmp}
      Performance profile of the total factorization floating-point operations for different scalings 
(all problems).
   }
    \vspace{3mm}

   \centering
   \input{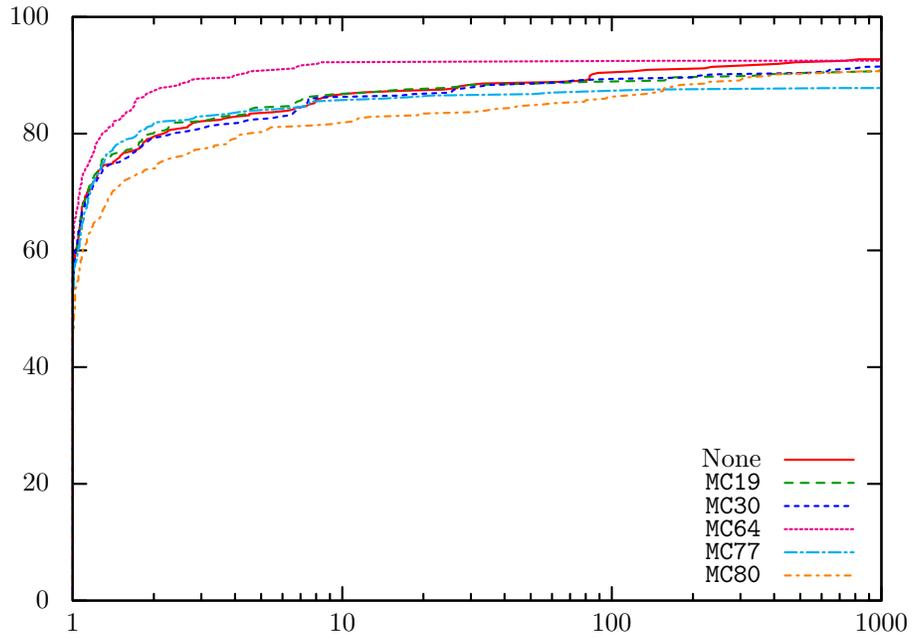}
\end{figure}

Closer examination of the results reveals that for each of the scalings there is little difference in
the number of Ipopt iterations used. For over 75\% of the problems, each of the scalings resulted in the
same number of iterations being used. For the majority of the remaining problems,
\texttt{HSL\_MC80} led to the smallest number of iterations, followed by \texttt{MC64}, with the
other scalings only a short distance behind. The main differences in the 
runtimes are explained by the number of delayed pivots that occur
(that is, the number of modifications that must be made to
the  pivot sequence selected by the analyse phase during the factorization
phase to ensure numerical stability)  and the number of floating-point operations
required to perform the factorizations. Figures~\ref{ndelay cmp} and \ref{nflops cmp} show performance
profiles of the average number of delayed pivots per factorization
and the total number of floating-point
operations, respectively. We observe that while \texttt{HSL\_MC80} is very
successful at eliminating delayed pivots, it does so at a substantial increase in the
number of floating-point operations required and this leads to it being
the most expensive option overall. However, \texttt{MC64} also
produces significantly fewer delayed pivots than the remaining scalings and is able to
leverage this to require fewer floating-point operations. The reason
      \texttt{MC64} is not the outright best scaling option  is that the call to \texttt{MC64}
      itself is computationally expensive. Furthermore, we have observed that, if the number
   of possible maximal matchings is small, the chosen scaling can be poor.

   \begin{table}[htbp]
      \caption{ \label{best scal probs}
         Ipopt runtimes (in seconds) for selected problems on which one scaling significantly outperforms 
   the others. A - indicates the run failed 
  (time or memory limit exceeded or failed to converge to a feasible solution). The numbers
   in brackets are the percentages of the runtime spent on scaling.
      }
       \vspace{3mm}

      \centering
      \small

      \begin{tabular}{l|*{6}{r@{ }l}}
         \hline
         Problem & \multicolumn{2}{c}{No scaling} & \multicolumn{2}{c}{\tt MC19} 
   & \multicolumn{2}{c}{\tt MC30} & \multicolumn{2}{c}{\tt MC64} 
   & \multicolumn{2}{c}{\tt MC77} & \multicolumn{2}{c}{\tt MC80} \\
         \hline
         SPIN & \bf 0.11&(0) & 0.14&(23) & 0.34&(12) & 0.19&(21) & 0.15&(13) & 0.20&(40) \\
         C-RELOAD & -& & \bf 0.95&(19) & -& & 2.13&(56) & -& & 2.18&(55) \\
         GRIDNETA & -& & 0.40&(17) & \bf 0.25&(22) & 0.43&(42) & 27.8&(29) & 0.44&(39) \\
         BRAINPC1 & 29.9&(0) & -& & -& & 44.2&(80) & \bf 5.88&(6) & -& \\
         CBRATU3D & 1.52&(0) & 1.50&(0.2) & 1.50&(0.2) & 1.52&(0.2) & 1.52&(0.2) & \bf 0.28&(5) \\
%         NCVXQP7 & 999&(0) & -& & -& & \bf 561&(21) & 665&(0.3) & -& \\
         STCQP1 & 30.5&(0) & 709&(0.4) & 8.47&(1) & \bf 3.68&(66) & -& & 7.80&(34) \\
         EIGENC & 161&(0) & 250&(0.2) & 238&(0.3) & 171&(4) & 317&(0.1) & \bf 68.3&(8) \\
         \hline
      \end{tabular}
   \end{table}

   \begin{table}[htbp]
      \caption{ \label{detail itr}
         Iteration counts for selected  problems using different scalings. Factorization
         counts are given in () brackets. A - indicates the run failed 
  (time or memory limit exceeded or failed to converge to a feasible solution).
      }
       \vspace{3mm}
      \centering
      \small
      \begin{tabular}{l|*{6}{r@{ }l}}
         \hline
         Problem & \multicolumn{2}{c}{No scaling} & \multicolumn{2}{c}{\tt MC19} 
   & \multicolumn{2}{c}{\tt MC30} & \multicolumn{2}{c}{\tt MC64} 
   & \multicolumn{2}{c}{\tt MC77} & \multicolumn{2}{c}{\tt MC80} \\
         \hline
         SPIN    & 7&(16)     & 7&(15)   & 9&(18)   & 7&(22)   & 7&(17)   &7&(14)\\
         C-RELOAD& -&         & 111&(242)& -&       & 116&(362)& -&       &123&(265)\\
         GRIDNETA& -&         & 27&(59)  & 20&(35)  & 21&(53)  & 1471&(4449)&20&(31) \\
         BRAINPC1& 82&(171)   & -&       & -&       & 60&(126) & 39&(77)  & -&   \\
         CBRATU3D& 3&(4)      & 3&(4)    & 3&(4)    &  3&(4)   & 3&(4)    & 3&(4)\\
%         NCVXQP7 & 158&(406)  & -&       & -&       & 129&(321)& 133&(343)& -&   \\
         STCQP1  & 14&(18)    & 123&(407)& 14&(18)  &  14&(18) & -&       &14&(14)\\
         EIGENC  & 14&(22)    & 14&(21)  & 14&(24)  &  13&(18) & 14&(20)  &13&(17)\\
         \hline
      \end{tabular}
   \end{table}

   \begin{table}[htbp]
      \caption{ \label{detail flops delayed pivots}
         Total factorization floating-point operation counts (in billions) for selected problems using
         different scalings. The largest number of delayed pivots (in thousands) for a single factorization is given in
         () brackets. A - indicates the run failed 
  (time or memory limit exceeded or failed to converge to a feasible solution).
      }
       \vspace{3mm}
      \centering
      \small
      \begin{tabular}{l|*{6}{r@{ }l}}
         \hline
         Problem & \multicolumn{2}{c}{No scaling} & \multicolumn{2}{c}{\tt MC19} & \multicolumn{2}{c}{\tt MC30} & \multicolumn{2}{c}{\tt MC64} & \multicolumn{2}{c}{\tt MC77} & \multicolumn{2}{c}{\tt MC80} \\
         \hline
         SPIN    & 0.049&(0) & 0.046&(0) & 1.03&(1) & 0.065&(0) & 0.052&(0) & 0.040&(0) \\
         C-RELOAD& -& & 0.660&(0.6) & -& & 0.933&(0.3) & -& & 0.684&(0.2) \\
         GRIDNETA& -& & 0.052&(0) & 0.031&(0) & 0.047&(0) & 4.11&(4) & 0.030&(0) \\
         BRAINPC1& 13.4&(756) & -& & -&          & 0.239&(0) & 0.151&(7) & -             \\
         CBRATU3D& 4.44&(7) & 4.44&(7) & 4.44&(7) & 4.44&(7) & 4.44&(7) & 0.670&(0) \\
%         NCVXQP7 & 2660&(147) & -&         & -& & 1220&(31) & 1620&(151) & -             \\
         STCQP1  & 178&(8) & 3390&(8) & 11.0&(8) & 1.49&(0) & -&       & 8.75&(0) \\
         EIGENC  & 945&(3) & 852&(3) & 1000&(3) & 746&(3) & 846&(3) & 387&(0) \\
         \hline
      \end{tabular}
   \end{table}

   \begin{table}
      \caption{ \label{detail scal time}
         Average scaling time (in milliseconds) per factorization.
   A - indicates the run failed 
  (time or memory limit exceeded or failed to converge to a feasible solution).
      }
       \vspace{3mm}
      \centering
      \small
      \begin{tabular}{l|*{5}{r}}
         \hline
         Problem & \tt MC19 & \tt MC30 & \tt MC64 & \tt MC77 & \tt MC80 \\
         \hline
         SPIN    & 2.2      & 2.3      & 2.4      & 1.1      & 5.7        \\
         C-RELOAD& 0.74     & -        & 3.3      & -        & 4.6        \\
         GRIDNETA& 1.2      & 1.6      & 3.4      & 1.8      & 5.4           \\
         BRAINPC1& 6.0      & -        & 281      & 4.6      & -             \\
         CBRATU3D& 0.87     & 1.2      & 1.0      & 1.1      & 3.6           \\
%         NCVXQP7 & -        & 69       & 370      & 5.8      & -             \\
         STCQP1  & 6.5      & 6.8      & 135      & -        & 148           \\
         EIGENC  & 28       & 32       & 343      & 18       & 326           \\
         \hline
      \end{tabular}
   \end{table}

   Tables~\ref{best scal probs}--\ref{detail scal time} report results
for a small selection of the test problems that were chosen to illustrate how
each of the scalings can outperform the others and to highlight the points of
the previous paragraph. Note that a single Ipopt iteration can result in
multiple factorizations as it attempts to determine the pertubation required to
achieve the correct inertia, or if iterative refinement fails and recovery
mechanisms are engaged.
For those problems with few delayed pivots and thus little difference in the
factorization cost per iteration (SPIN, C-RELOAD, GRIDNETA, BRAINPC1),
the scaling that achieves the minimum number of iterations with the smallest
scaling time overhead is the fastest. For those problems where delayed pivots are
significant (CBRATU3D, STCQP1, EIGNEC), the more expensive \texttt{MC64} and
\texttt{HSL\_MC80} scalings that reduce the number of delayed pivots lead to
the best runtimes.

These findings may be compared with previous research by Hogg and Scott
\cite{hosc:2008a,hosc:2012b} into the effects of scalings on a large set of sparse
linear systems arising from a wide range of practical applications. They 
found that for tough indefinite systems \texttt{HSL\_MC80} produces
the highest quality scalings, but takes significantly longer than
non-matching based scalings to run. In some cases, the time for \texttt{HSL\_MC80} and \texttt{MC64} dominated
the cost of the linear solver. Hogg and Scott also reported that  for these
systems, \texttt{MC30}
is generally not competitive with \texttt{MC77} in terms of  quality and run-time.

\section{Effects of dynamic scaling}\label{dynamic scaling}

\begin{figure}[htbp]
   \caption{ \label{on demand cmp}
      Performance profile of times for different scalings used 
with the on demand heuristic (problems that enable scaling only).
   }
   \centering
   \begin{tikzpicture}[gnuplot]
%% generated with GNUPLOT 4.4p3 (Lua 5.1.5; terminal rev. 97, script rev. 96a)
%% Mon 10 Dec 2012 12:35:30 GMT
\gpcolor{gp lt color border}
\gpsetlinetype{gp lt border}
\gpsetlinewidth{2.00}
\draw[gp path] (1.196,0.616)--(1.376,0.616);
\draw[gp path] (11.947,0.616)--(11.767,0.616);
\node[gp node right] at (1.012,0.616) { 0};
\draw[gp path] (1.196,2.169)--(1.376,2.169);
\draw[gp path] (11.947,2.169)--(11.767,2.169);
\node[gp node right] at (1.012,2.169) { 20};
\draw[gp path] (1.196,3.722)--(1.376,3.722);
\draw[gp path] (11.947,3.722)--(11.767,3.722);
\node[gp node right] at (1.012,3.722) { 40};
\draw[gp path] (1.196,5.275)--(1.376,5.275);
\draw[gp path] (11.947,5.275)--(11.767,5.275);
\node[gp node right] at (1.012,5.275) { 60};
\draw[gp path] (1.196,6.828)--(1.376,6.828);
\draw[gp path] (11.947,6.828)--(11.767,6.828);
\node[gp node right] at (1.012,6.828) { 80};
\draw[gp path] (1.196,8.381)--(1.376,8.381);
\draw[gp path] (11.947,8.381)--(11.767,8.381);
\node[gp node right] at (1.012,8.381) { 100};
\draw[gp path] (1.196,0.616)--(1.196,0.796);
\draw[gp path] (1.196,8.381)--(1.196,8.201);
\node[gp node center] at (1.196,0.308) { 1};
\draw[gp path] (3.884,0.616)--(3.884,0.796);
\draw[gp path] (3.884,8.381)--(3.884,8.201);
\node[gp node center] at (3.884,0.308) { 1.5};
\draw[gp path] (6.572,0.616)--(6.572,0.796);
\draw[gp path] (6.572,8.381)--(6.572,8.201);
\node[gp node center] at (6.572,0.308) { 2};
\draw[gp path] (9.259,0.616)--(9.259,0.796);
\draw[gp path] (9.259,8.381)--(9.259,8.201);
\node[gp node center] at (9.259,0.308) { 2.5};
\draw[gp path] (11.947,0.616)--(11.947,0.796);
\draw[gp path] (11.947,8.381)--(11.947,8.201);
\node[gp node center] at (11.947,0.308) { 3};
\draw[gp path] (1.196,8.381)--(1.196,0.616)--(11.947,0.616)--(11.947,8.381)--cycle;
\node[gp node right] at (10.479,2.490) {None};
\gpcolor{gp lt color 0}
\gpsetlinetype{gp lt plot 0}
\draw[gp path] (10.663,2.490)--(11.579,2.490);
\draw[gp path] (1.196,0.741)--(1.196,0.867)--(1.196,0.992)--(1.196,1.117)--(1.196,1.242)%
  --(1.196,1.368)--(1.196,1.493)--(1.250,1.618)--(1.626,1.743)--(1.895,1.868)--(2.002,1.994)%
  --(2.110,2.119)--(2.486,2.244)--(2.916,2.369)--(2.916,2.494)--(3.239,2.620)--(3.507,2.745)%
  --(3.669,2.870)--(3.776,2.996)--(3.884,3.121)--(3.938,3.246)--(3.938,3.371)--(4.421,3.497)%
  --(4.636,3.622)--(4.905,3.747)--(5.281,3.873)--(5.873,3.998)--(7.109,4.123)--(7.163,4.248)%
  --(8.292,4.373)--(11.839,4.499)--(11.947,4.500);
\gpcolor{gp lt color border}
\node[gp node right] at (10.479,2.182) {\tt MC19-OD};
\gpcolor{gp lt color 1}
\gpsetlinetype{gp lt plot 1}
\draw[gp path] (10.663,2.182)--(11.579,2.182);
\draw[gp path] (1.196,0.741)--(1.196,0.867)--(1.196,0.992)--(1.196,1.117)--(1.196,1.242)%
  --(1.196,1.368)--(1.196,1.493)--(1.196,1.618)--(1.250,1.743)--(1.304,1.868)--(1.357,1.994)%
  --(1.626,2.119)--(1.680,2.244)--(1.680,2.369)--(1.841,2.494)--(2.164,2.620)--(2.217,2.745)%
  --(2.916,2.870)--(2.916,2.996)--(3.292,3.121)--(3.669,3.246)--(3.938,3.371)--(5.389,3.497)%
  --(5.604,3.622)--(6.088,3.747)--(7.647,3.873)--(7.754,3.998)--(7.915,4.123)--(11.947,4.212);
\gpcolor{gp lt color border}
\node[gp node right] at (10.479,1.874) {\tt MC30-OD};
\gpcolor{gp lt color 2}
\gpsetlinetype{gp lt plot 2}
\draw[gp path] (10.663,1.874)--(11.579,1.874);
\draw[gp path] (1.196,0.741)--(1.196,0.867)--(1.196,0.992)--(1.304,1.117)--(1.411,1.242)%
  --(1.519,1.368)--(1.572,1.493)--(1.572,1.618)--(1.680,1.743)--(1.787,1.868)--(1.841,1.994)%
  --(2.110,2.119)--(2.379,2.244)--(2.647,2.369)--(3.077,2.494)--(3.239,2.620)--(3.292,2.745)%
  --(3.346,2.870)--(3.830,2.996)--(3.830,3.121)--(5.281,3.246)--(5.443,3.371)--(5.658,3.497)%
  --(6.141,3.622)--(7.163,3.747)--(9.044,3.873)--(9.743,3.998)--(10.119,4.123)--(10.119,4.248)%
  --(10.119,4.373)--(10.173,4.499)--(10.388,4.624)--(11.947,4.653);
\gpcolor{gp lt color border}
\node[gp node right] at (10.479,1.566) {\tt MC64-OD};
\gpcolor{gp lt color 3}
\gpsetlinetype{gp lt plot 3}
\draw[gp path] (10.663,1.566)--(11.579,1.566);
\draw[gp path] (1.196,0.741)--(1.196,0.867)--(1.196,0.992)--(1.196,1.117)--(1.196,1.242)%
  --(1.196,1.368)--(1.196,1.493)--(1.196,1.618)--(1.196,1.743)--(1.196,1.868)--(1.196,1.994)%
  --(1.196,2.119)--(1.304,2.244)--(1.411,2.369)--(1.465,2.494)--(1.519,2.620)--(1.519,2.745)%
  --(1.572,2.870)--(1.626,2.996)--(1.680,3.121)--(1.734,3.246)--(2.110,3.371)--(2.164,3.497)%
  --(2.432,3.622)--(4.045,3.747)--(4.744,3.873)--(5.228,3.998)--(5.604,4.123)--(7.485,4.248)%
  --(7.593,4.373)--(10.764,4.499)--(11.947,4.552);
\gpcolor{gp lt color border}
\node[gp node right] at (10.479,1.258) {\tt MC77-OD};
\gpcolor{gp lt color 4}
\gpsetlinetype{gp lt plot 4}
\draw[gp path] (10.663,1.258)--(11.579,1.258);
\draw[gp path] (1.196,0.741)--(1.196,0.867)--(1.196,0.992)--(1.196,1.117)--(1.196,1.242)%
  --(1.196,1.368)--(1.196,1.493)--(1.196,1.618)--(1.196,1.743)--(1.304,1.868)--(1.465,1.994)%
  --(1.519,2.119)--(1.519,2.244)--(1.626,2.369)--(1.680,2.494)--(1.895,2.620)--(2.164,2.745)%
  --(2.325,2.870)--(2.647,2.996)--(3.024,3.121)--(3.131,3.246)--(3.185,3.371)--(3.346,3.497)%
  --(5.550,3.622)--(5.711,3.747)--(6.787,3.873)--(7.539,3.998)--(8.614,4.123)--(8.829,4.248)%
  --(11.087,4.373)--(11.947,4.386);
\gpcolor{gp lt color border}
\node[gp node right] at (10.479,0.950) {\tt MC80-OD};
\gpcolor{gp lt color 5}
\gpsetlinetype{gp lt plot 5}
\draw[gp path] (10.663,0.950)--(11.579,0.950);
\draw[gp path] (1.196,0.741)--(1.196,0.867)--(1.196,0.992)--(1.196,1.117)--(1.196,1.242)%
  --(1.196,1.368)--(1.196,1.493)--(1.196,1.618)--(1.196,1.743)--(1.196,1.868)--(1.196,1.994)%
  --(1.465,2.119)--(1.519,2.244)--(1.626,2.369)--(1.626,2.494)--(1.734,2.620)--(1.787,2.745)%
  --(2.002,2.870)--(2.432,2.996)--(2.540,3.121)--(2.755,3.246)--(2.809,3.371)--(3.131,3.497)%
  --(3.185,3.622)--(4.636,3.747)--(5.066,3.873)--(5.120,3.998)--(5.819,4.123)--(9.313,4.248)%
  --(10.979,4.373)--(11.409,4.499)--(11.678,4.624)--(11.947,4.629);
\gpcolor{gp lt color border}
\gpsetlinetype{gp lt border}
\draw[gp path] (1.196,8.381)--(1.196,0.616)--(11.947,0.616)--(11.947,8.381)--cycle;
%% coordinates of the plot area
\gpdefrectangularnode{gp plot 1}{\pgfpoint{1.196cm}{0.616cm}}{\pgfpoint{11.947cm}{8.381cm}}
\end{tikzpicture}
\end{figure}
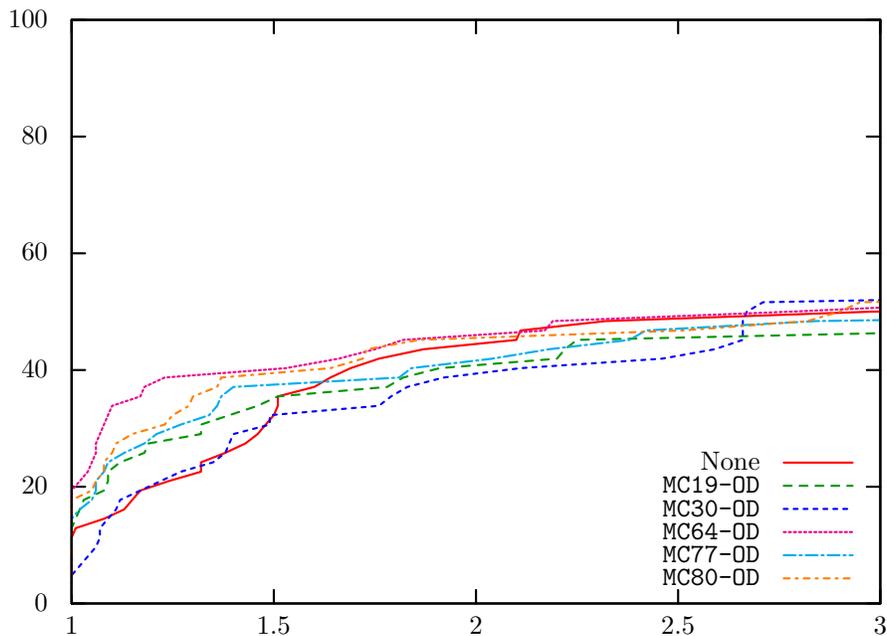

Due to the high overhead of scaling (and particularly of the matching-based scalings), it is logical to experiment
with the dynamic use and reuse of scalings. Currently, the default scaling strategy
when using an HSL linear solver within Ipopt is to start the computation without scaling 
and to switch to using scaling when iterative refinement fails to converge to
the required accuracy.
%\footnote{Note that, by default, the HSL solver {\tt MA57}
%uses {\tt MC64} scaling and, if scaling within Ipopt is not required, the scaling control 
%for {\tt MA57} needs to be reset.}.
At this point, to improve the accuracy of subsequent factorizations,
the pivot threshold parameter $u$ that is used
by the linear solver to control numerical stability
while limiting the fill in the computed factors, is increased 
and all subsequent
factorizations use \texttt{MC19} scaling. In keeping with Ipopt's
configuration parameters, we refer to this as the ``on demand'' heuristic.
Examination of our test set shows that, using {\tt HSL\_MA97},
scaling is switched on in this way 
for only 62 out of the 386 problems (so that, for the majority of
the test problems, scaling is not used). This confirms our earlier
finding that for the entire test set no scaling  is
the winning strategy (Figure~\ref{full scaling cmp}).
For the subset of 62 problems, Figure~\ref{on demand cmp} compares the results of
using on demand scaling with each of our scalings. 
When no scaling is used, the failed iterative refinement is resolved by a
perturbation of the optimization parameters.
We observe that we are  able to solve around half of these tough problems.

We have established that the effect of scaling is stronger than merely
increased accuracy: it can also be used to reduce the number of delayed pivots
(and thus the fill in the factors and the flop count).
This suggests the introduction of a second heuristic related to this.
Further, it is not clear whether it is necessary to recompute the scaling at
\textit{each} iteration, or if we can reuse the scaling computed at an
earlier iteration.
We therefore modified our {\tt HSL\_MA97} Ipopt driver to allow a number of different 
heuristics to be tested.
We use the following notation to
denote these:
\begin{description}
   \item[{\tt MCXX-OD}] ``on demand'' heuristic, where {\tt MCXX} scaling is computed
      for subsequent factorizations after the first failure of iterative refinement.
   \item[{\tt MCXX-ODR}] ``on demand reuse'' heuristic, as {\tt MCXX-OD}, however the
      scaling is only computed for the next factorization, and then reused
      thereafter (until iterative refinement again fails to converge).
   \item[{\tt MCXX-HD}] ``high delay'' heuristic, where {\tt MCXX} scaling is computed
      for all subsequent factorization after the first factorization for which
      more than $0.05n$ delayed pivots occur.
   \item[{\tt MCXX-HDR}] ``high delay reuse'' heuristic, as {\tt MCXX-HD}, however the
      scaling is only computed on the next factorization, and then reused
      thereafter. The number of delayed pivots reported by that next factorization is
      recorded and, if for a subsequent factorization more than $0.05n$ 
      {\it additional} delayed pivots are encountered, a new scaling is computed for
      the following factorization.
   \item[{\tt MCXX-ODHD}] ``on demand/high delay'' heuristic, where {\tt MCXX} scaling
      is enabled for all subsequent factorizations if either of the {\tt MCXX-OD} or
      {\tt MCXX-HD} conditions are met.
   \item[{\tt MCXX-ODHDR}] ``on demand/high delay reuse'' heuristic, where {\tt MCXX}
      scaling is enabled for some subsequent factorizations based on a
      combination of the {\tt MCXX-ODR} and {\tt MCXX-HDR} conditions.
\end{description}

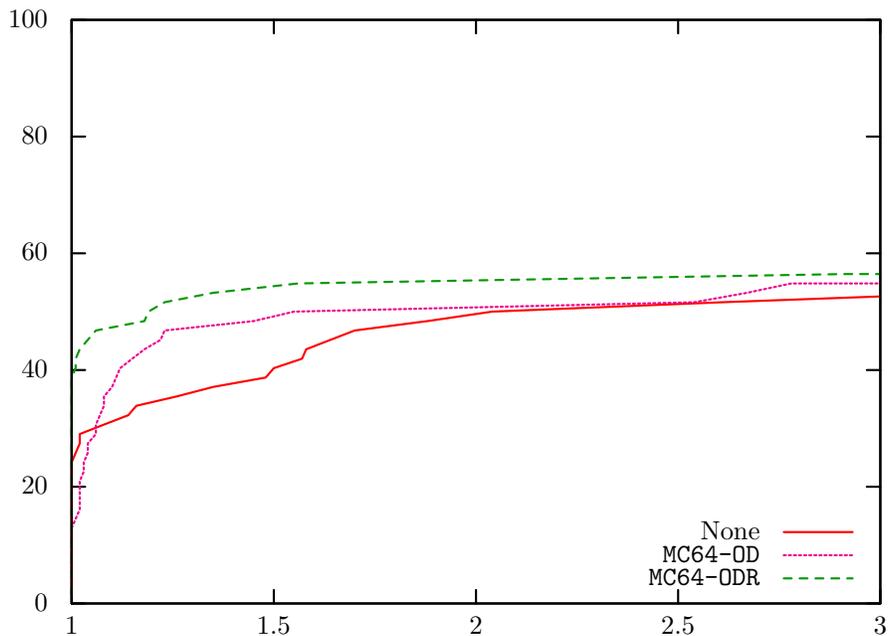
\begin{figure}
   \caption{ \label{fig:od vs odr}
      Performance profile of times for different {\tt MC64}-based heuristics (problems that enable scaling only).
   }
   \centering
   \begin{tikzpicture}[gnuplot]
%% generated with GNUPLOT 4.4p3 (Lua 5.1.5; terminal rev. 97, script rev. 96a)
%% Mon 10 Dec 2012 12:45:27 GMT
\gpcolor{gp lt color border}
\gpsetlinetype{gp lt border}
\gpsetlinewidth{2.00}
\draw[gp path] (1.196,0.616)--(1.376,0.616);
\draw[gp path] (11.947,0.616)--(11.767,0.616);
\node[gp node right] at (1.012,0.616) { 0};
\draw[gp path] (1.196,2.169)--(1.376,2.169);
\draw[gp path] (11.947,2.169)--(11.767,2.169);
\node[gp node right] at (1.012,2.169) { 20};
\draw[gp path] (1.196,3.722)--(1.376,3.722);
\draw[gp path] (11.947,3.722)--(11.767,3.722);
\node[gp node right] at (1.012,3.722) { 40};
\draw[gp path] (1.196,5.275)--(1.376,5.275);
\draw[gp path] (11.947,5.275)--(11.767,5.275);
\node[gp node right] at (1.012,5.275) { 60};
\draw[gp path] (1.196,6.828)--(1.376,6.828);
\draw[gp path] (11.947,6.828)--(11.767,6.828);
\node[gp node right] at (1.012,6.828) { 80};
\draw[gp path] (1.196,8.381)--(1.376,8.381);
\draw[gp path] (11.947,8.381)--(11.767,8.381);
\node[gp node right] at (1.012,8.381) { 100};
\draw[gp path] (1.196,0.616)--(1.196,0.796);
\draw[gp path] (1.196,8.381)--(1.196,8.201);
\node[gp node center] at (1.196,0.308) { 1};
\draw[gp path] (3.884,0.616)--(3.884,0.796);
\draw[gp path] (3.884,8.381)--(3.884,8.201);
\node[gp node center] at (3.884,0.308) { 1.5};
\draw[gp path] (6.572,0.616)--(6.572,0.796);
\draw[gp path] (6.572,8.381)--(6.572,8.201);
\node[gp node center] at (6.572,0.308) { 2};
\draw[gp path] (9.259,0.616)--(9.259,0.796);
\draw[gp path] (9.259,8.381)--(9.259,8.201);
\node[gp node center] at (9.259,0.308) { 2.5};
\draw[gp path] (11.947,0.616)--(11.947,0.796);
\draw[gp path] (11.947,8.381)--(11.947,8.201);
\node[gp node center] at (11.947,0.308) { 3};
\draw[gp path] (1.196,8.381)--(1.196,0.616)--(11.947,0.616)--(11.947,8.381)--cycle;
\node[gp node right] at (10.479,1.566) {None};
\gpcolor{gp lt color 0}
\gpsetlinetype{gp lt plot 0}
\draw[gp path] (10.663,1.566)--(11.579,1.566);
\draw[gp path] (1.196,0.741)--(1.196,0.867)--(1.196,0.992)--(1.196,1.117)--(1.196,1.242)%
  --(1.196,1.368)--(1.196,1.493)--(1.196,1.618)--(1.196,1.743)--(1.196,1.868)--(1.196,1.994)%
  --(1.196,2.119)--(1.196,2.244)--(1.196,2.369)--(1.196,2.494)--(1.250,2.620)--(1.304,2.745)%
  --(1.304,2.870)--(1.626,2.996)--(1.949,3.121)--(2.056,3.246)--(2.594,3.371)--(3.077,3.497)%
  --(3.776,3.622)--(3.884,3.747)--(4.260,3.873)--(4.314,3.998)--(4.636,4.123)--(4.959,4.248)%
  --(5.926,4.373)--(6.787,4.499)--(9.851,4.624)--(11.947,4.701);
\gpcolor{gp lt color border}
\node[gp node right] at (10.479,1.258) {\tt MC64-OD};
\gpcolor{gp lt color 3}
\gpsetlinetype{gp lt plot 3}
\draw[gp path] (10.663,1.258)--(11.579,1.258);
\draw[gp path] (1.196,0.741)--(1.196,0.867)--(1.196,0.992)--(1.196,1.117)--(1.196,1.242)%
  --(1.196,1.368)--(1.196,1.493)--(1.196,1.618)--(1.250,1.743)--(1.304,1.868)--(1.304,1.994)%
  --(1.304,2.119)--(1.304,2.244)--(1.357,2.369)--(1.357,2.494)--(1.411,2.620)--(1.411,2.745)%
  --(1.519,2.870)--(1.519,2.996)--(1.572,3.121)--(1.626,3.246)--(1.626,3.371)--(1.734,3.497)%
  --(1.787,3.622)--(1.841,3.747)--(2.002,3.873)--(2.164,3.998)--(2.379,4.123)--(2.432,4.248)%
  --(3.615,4.373)--(4.153,4.499)--(9.474,4.624)--(10.173,4.749)--(10.764,4.874)--(11.947,4.874);
\gpcolor{gp lt color border}
\node[gp node right] at (10.479,0.950) {\tt MC64-ODR};
\gpcolor{gp lt color 1}
\gpsetlinetype{gp lt plot 1}
\draw[gp path] (10.663,0.950)--(11.579,0.950);
\draw[gp path] (1.196,0.741)--(1.196,0.867)--(1.196,0.992)--(1.196,1.117)--(1.196,1.242)%
  --(1.196,1.368)--(1.196,1.493)--(1.196,1.618)--(1.196,1.743)--(1.196,1.868)--(1.196,1.994)%
  --(1.196,2.119)--(1.196,2.244)--(1.196,2.369)--(1.196,2.494)--(1.196,2.620)--(1.196,2.745)%
  --(1.196,2.870)--(1.196,2.996)--(1.196,3.121)--(1.196,3.246)--(1.196,3.371)--(1.196,3.497)%
  --(1.196,3.622)--(1.250,3.747)--(1.250,3.873)--(1.304,3.998)--(1.411,4.123)--(1.519,4.248)%
  --(2.164,4.373)--(2.217,4.499)--(2.432,4.624)--(3.077,4.749)--(4.206,4.874)--(11.463,4.999)%
  --(11.947,4.999);
\gpcolor{gp lt color border}
\gpsetlinetype{gp lt border}
\draw[gp path] (1.196,8.381)--(1.196,0.616)--(11.947,0.616)--(11.947,8.381)--cycle;
%% coordinates of the plot area
\gpdefrectangularnode{gp plot 1}{\pgfpoint{1.196cm}{0.616cm}}{\pgfpoint{11.947cm}{8.381cm}}
\end{tikzpicture}
\end{figure}

For the on demand heuristic with {\tt MC64} scaling, Figure~\ref{fig:od vs odr} 
illustrates the benefits of
reusing the scaling (again for the subset of
62 tough problems). For the large problems in  this subset we get some worthwhile savings
resulting from the reduction in the total cost of scaling
without a deterioration in the factorization quality.
For each heuristic, similar behaviour was observed with the reuse version outperforming
the  version that recomputes the scaling. For the complete test set and
the subset of large problems
Figures~\ref{fig:mc64 heuristic cmp} and \ref{fig:mc64 heuristic cmp big}, respectively,
compare the different reuse heuristics; reliability results are given in Table~\ref{reliability final}.
These show that using scaling to alleviate the number of delayed pivots
significantly increases reliability (by preventing the run aborting
 because the time or memory limits are exceeded), while the conditions based on
failure of iterative refinement also assist with reliability but are less effective. These latter
conditions also result in a moderate performance advantage for a number of problems.
The optimal strategy for the larger problems appears to be 
{\tt MC64-ODHDR}, while for the entire test set (dominated by small problems), the best heuristic
is {\tt MC64-ODR}. This is because the overhead of computing the
scaling on the small problems is greater than the saving in the factorization time
resulting from reducing the number of delayed pivots.
%\begin{figure}
%   \caption{ \label{fig:od vs odr big}
%      Performance profile of times for different {\tt MC64}-based heuristics (large problems).
%   }
%   \centering
%   \input{mc64_od_vs_odr_big.tex}
%\end{figure}
\begin{figure}
   \caption{ \label{fig:mc64 heuristic cmp}
      Performance profile of times for different {\tt MC64}-based heuristics (all problems).
   }
   \centering
   \input{heuristic.tex}
\end{figure}
\begin{figure}
   \caption{ \label{fig:mc64 heuristic cmp big}
      Performance profile of times for different {\tt MC64}-based heuristics (large problems).
   }
   \centering
   \begin{tikzpicture}[gnuplot]
%% generated with GNUPLOT 4.4p3 (Lua 5.1.5; terminal rev. 97, script rev. 96a)
%% Mon 10 Dec 2012 13:23:41 GMT
\gpcolor{gp lt color border}
\gpsetlinetype{gp lt border}
\gpsetlinewidth{2.00}
\draw[gp path] (1.196,0.616)--(1.376,0.616);
\draw[gp path] (11.947,0.616)--(11.767,0.616);
\node[gp node right] at (1.012,0.616) { 0};
\draw[gp path] (1.196,2.169)--(1.376,2.169);
\draw[gp path] (11.947,2.169)--(11.767,2.169);
\node[gp node right] at (1.012,2.169) { 20};
\draw[gp path] (1.196,3.722)--(1.376,3.722);
\draw[gp path] (11.947,3.722)--(11.767,3.722);
\node[gp node right] at (1.012,3.722) { 40};
\draw[gp path] (1.196,5.275)--(1.376,5.275);
\draw[gp path] (11.947,5.275)--(11.767,5.275);
\node[gp node right] at (1.012,5.275) { 60};
\draw[gp path] (1.196,6.828)--(1.376,6.828);
\draw[gp path] (11.947,6.828)--(11.767,6.828);
\node[gp node right] at (1.012,6.828) { 80};
\draw[gp path] (1.196,8.381)--(1.376,8.381);
\draw[gp path] (11.947,8.381)--(11.767,8.381);
\node[gp node right] at (1.012,8.381) { 100};
\draw[gp path] (1.196,0.616)--(1.196,0.796);
\draw[gp path] (1.196,8.381)--(1.196,8.201);
\node[gp node center] at (1.196,0.308) { 1};
\draw[gp path] (3.884,0.616)--(3.884,0.796);
\draw[gp path] (3.884,8.381)--(3.884,8.201);
\node[gp node center] at (3.884,0.308) { 1.5};
\draw[gp path] (6.572,0.616)--(6.572,0.796);
\draw[gp path] (6.572,8.381)--(6.572,8.201);
\node[gp node center] at (6.572,0.308) { 2};
\draw[gp path] (9.259,0.616)--(9.259,0.796);
\draw[gp path] (9.259,8.381)--(9.259,8.201);
\node[gp node center] at (9.259,0.308) { 2.5};
\draw[gp path] (11.947,0.616)--(11.947,0.796);
\draw[gp path] (11.947,8.381)--(11.947,8.201);
\node[gp node center] at (11.947,0.308) { 3};
\draw[gp path] (1.196,8.381)--(1.196,0.616)--(11.947,0.616)--(11.947,8.381)--cycle;
\node[gp node right] at (10.479,1.874) {None};
\gpcolor{gp lt color 0}
\gpsetlinetype{gp lt plot 0}
\draw[gp path] (10.663,1.874)--(11.579,1.874);
\draw[gp path] (1.196,0.760)--(1.196,0.903)--(1.196,1.048)--(1.196,1.191)--(1.196,1.335)%
  --(1.196,1.479)--(1.196,1.622)--(1.196,1.766)--(1.196,1.910)--(1.196,2.054)--(1.196,2.198)%
  --(1.196,2.341)--(1.196,2.485)--(1.196,2.629)--(1.196,2.773)--(1.196,2.917)--(1.196,3.060)%
  --(1.196,3.204)--(1.196,3.349)--(1.196,3.492)--(1.196,3.636)--(1.196,3.779)--(1.196,3.923)%
  --(1.196,4.067)--(1.196,4.211)--(1.196,4.355)--(1.196,4.499)--(1.196,4.642)--(1.250,4.786)%
  --(1.250,4.930)--(1.250,5.074)--(1.250,5.218)--(1.250,5.361)--(1.250,5.505)--(1.250,5.648)%
  --(1.250,5.793)--(1.250,5.937)--(1.304,6.080)--(1.357,6.224)--(1.357,6.368)--(1.411,6.512)%
  --(1.411,6.656)--(1.519,6.799)--(1.626,6.943)--(2.271,7.087)--(2.540,7.231)--(2.647,7.375)%
  --(5.926,7.518)--(11.947,7.649);
\gpcolor{gp lt color border}
\node[gp node right] at (10.479,1.566) {\tt MC64-ODR};
\gpcolor{gp lt color 1}
\gpsetlinetype{gp lt plot 1}
\draw[gp path] (10.663,1.566)--(11.579,1.566);
\draw[gp path] (1.196,0.760)--(1.196,0.903)--(1.196,1.048)--(1.196,1.191)--(1.196,1.335)%
  --(1.196,1.479)--(1.196,1.622)--(1.196,1.766)--(1.196,1.910)--(1.196,2.054)--(1.196,2.198)%
  --(1.196,2.341)--(1.196,2.485)--(1.196,2.629)--(1.196,2.773)--(1.196,2.917)--(1.196,3.060)%
  --(1.196,3.204)--(1.196,3.349)--(1.196,3.492)--(1.196,3.636)--(1.196,3.779)--(1.196,3.923)%
  --(1.196,4.067)--(1.196,4.211)--(1.196,4.355)--(1.250,4.499)--(1.250,4.642)--(1.250,4.786)%
  --(1.250,4.930)--(1.250,5.074)--(1.304,5.218)--(1.304,5.361)--(1.304,5.505)--(1.304,5.648)%
  --(1.304,5.793)--(1.357,5.937)--(1.357,6.080)--(1.357,6.224)--(1.357,6.368)--(1.411,6.512)%
  --(1.411,6.656)--(1.411,6.799)--(1.519,6.943)--(1.626,7.087)--(1.680,7.231)--(1.734,7.375)%
  --(2.002,7.518)--(2.647,7.662)--(11.947,7.798);
\gpcolor{gp lt color border}
\node[gp node right] at (10.479,1.258) {\tt MC64-HDR};
\gpcolor{gp lt color 2}
\gpsetlinetype{gp lt plot 2}
\draw[gp path] (10.663,1.258)--(11.579,1.258);
\draw[gp path] (1.196,0.760)--(1.196,0.903)--(1.196,1.048)--(1.196,1.191)--(1.196,1.335)%
  --(1.196,1.479)--(1.196,1.622)--(1.196,1.766)--(1.196,1.910)--(1.196,2.054)--(1.196,2.198)%
  --(1.196,2.341)--(1.196,2.485)--(1.196,2.629)--(1.196,2.773)--(1.196,2.917)--(1.196,3.060)%
  --(1.196,3.204)--(1.196,3.349)--(1.196,3.492)--(1.196,3.636)--(1.196,3.779)--(1.196,3.923)%
  --(1.250,4.067)--(1.250,4.211)--(1.250,4.355)--(1.250,4.499)--(1.250,4.642)--(1.250,4.786)%
  --(1.250,4.930)--(1.250,5.074)--(1.250,5.218)--(1.250,5.361)--(1.250,5.505)--(1.304,5.648)%
  --(1.304,5.793)--(1.304,5.937)--(1.304,6.080)--(1.304,6.224)--(1.304,6.368)--(1.357,6.512)%
  --(1.357,6.656)--(1.411,6.799)--(1.519,6.943)--(1.626,7.087)--(1.734,7.231)--(1.734,7.375)%
  --(1.734,7.518)--(1.787,7.662)--(2.432,7.806)--(2.486,7.949)--(2.540,8.094)--(11.947,8.094);
\gpcolor{gp lt color border}
\node[gp node right] at (10.479,0.950) {\tt MC64-ODHDR};
\gpcolor{gp lt color 3}
\gpsetlinetype{gp lt plot 3}
\draw[gp path] (10.663,0.950)--(11.579,0.950);
\draw[gp path] (1.196,0.760)--(1.196,0.903)--(1.196,1.048)--(1.196,1.191)--(1.196,1.335)%
  --(1.196,1.479)--(1.196,1.622)--(1.196,1.766)--(1.196,1.910)--(1.196,2.054)--(1.196,2.198)%
  --(1.196,2.341)--(1.196,2.485)--(1.196,2.629)--(1.196,2.773)--(1.196,2.917)--(1.196,3.060)%
  --(1.196,3.204)--(1.196,3.349)--(1.196,3.492)--(1.196,3.636)--(1.196,3.779)--(1.196,3.923)%
  --(1.196,4.067)--(1.196,4.211)--(1.196,4.355)--(1.196,4.499)--(1.196,4.642)--(1.196,4.786)%
  --(1.250,4.930)--(1.250,5.074)--(1.250,5.218)--(1.250,5.361)--(1.250,5.505)--(1.250,5.648)%
  --(1.250,5.793)--(1.250,5.937)--(1.250,6.080)--(1.250,6.224)--(1.304,6.368)--(1.304,6.512)%
  --(1.304,6.656)--(1.304,6.799)--(1.304,6.943)--(1.357,7.087)--(1.411,7.231)--(1.465,7.375)%
  --(1.519,7.518)--(1.734,7.662)--(2.379,7.806)--(2.486,7.949)--(2.486,8.094)--(11.947,8.094);
\gpcolor{gp lt color border}
\gpsetlinetype{gp lt border}
\draw[gp path] (1.196,8.381)--(1.196,0.616)--(11.947,0.616)--(11.947,8.381)--cycle;
%% coordinates of the plot area
\gpdefrectangularnode{gp plot 1}{\pgfpoint{1.196cm}{0.616cm}}{\pgfpoint{11.947cm}{8.381cm}}
\end{tikzpicture}
%% gnuplot variables
\end{figure}
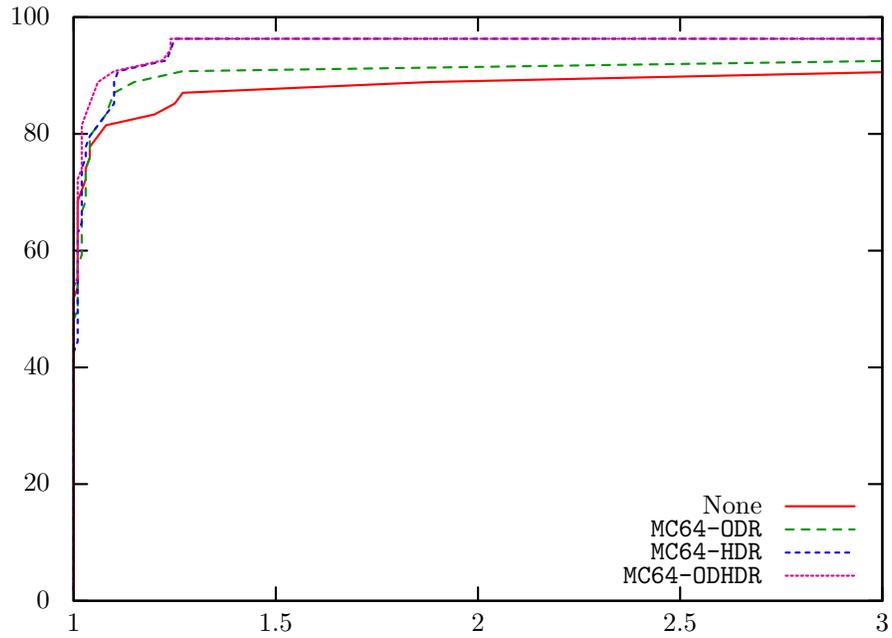
\begin{table}
   \caption{ \label{reliability final}
      Reliability of different scalings used on the set of all (non-trivial, tractable) problems and the subset of large problems.
   }
    \vspace{3mm}

   \centering
   \small

   \begin{tabular}{l|cc}
      \hline
                        & All    & Large \\
                        & problems & problems \\
      \hline
      None              & 92.8\% & 94.4\% \\
      \texttt{MC64}     & 92.5\% & 96.3\% \\
      \texttt{MC64-OD}  & 92.2\% & 94.4\% \\ 
      \texttt{MC64-ODR} & 92.5\% & 94.4\% \\
      \texttt{MC64-HDR} & 93.0\% & 96.3\% \\
      \texttt{MC64-ODHDR}&93.5\% & 96.3\% \\
      \hline
   \end{tabular}
\end{table}

\section{Choice of scaling for specific problem(s)} \label{advice}
We recommend that the initial choice of scaling be the {\tt MC64-ODHDR} method; this is
the default setting for our {\tt HSL\_MA97} driver within Ipopt. If the performance is not satisfactory,
the user can enable printing of factorization statistics and information on the use of scaling by selecting
the appropriate Ipopt output options.

If there are a large number of delays on the first iteration, an option exists in our driver
to force the use of {\tt MC64} on the first iteration in addition to the default heuristic.
If a large number of delays occur on iterations that are reusing the scaling computed on
an earlier iteration, the user is advised to investigate performance using the {\tt MC64-ODHD} heuristic.
Should a large number of delays be present on an iteration where the {\tt MC64} scaling has been
computed, use of {\tt HSL\_MC80} scaling should be considered by enabling the appropriate ordering option.

We have observed that, in some cases where there is a small number of possible maximum matchings across which
{\tt MC64} can optimize, the computed scaling factors can be poor. If this behaviour is suspected (for example,
if a warning has been issued), use of alternative scalings (such as {\tt MC77}) should be tested.
Alternative scalings (or no scaling) may also be considered for small problems for which computing
the matching-based scalings can add an unacceptably large overhead.

For cases where top performance is required, we recommend comprehensive testing using all possible scalings
(including no scaling) and heuristics that are included in our driver 
so as to determine the best option for the specific class
of problems involved.

\section{Conclusions} \label{conclusions}
In this study, we have used the CUTEr test set to explore the effects of scalings
within Ipopt. This set is dominated by small problems. For these, the factorization
time per iteration is a small proportion of the total Ipopt runtime and so
the potential gains from scaling are limited. In particular, an expensive matching-based scaling,
while improving reliability, can dominate the total runtime.
For larger problems within the CUTEr set, we have seen that the use of scaling can offer worthwhile savings
but this is highly problem dependent. Heuristics based, for example, on the time taken to compute the
scaling on the first iteration could be considered but time-based heuristics are likely
to be very machine dependent. Instead, we have proposed a number of heuristics based on enabling scaling when either
iterative refinement fails or the number of delayed pivots becomes large. To reduce
scaling costs, we have also considered the reuse of scaling and showed that this
can lead to reductions in the runtime. We would welcome the opportunity to run our tests
on other large-scale problems, particularly tough problems that are not necessarily well scaled;
we are always seeking to expand the set of test problems that we have available.

All our reported timings were for runs performed in serial. {\tt HSL\_MA97} is, however, a parallel
solver. Currently, only serial implementations of the considered scalings are available.
If we run Ipopt with  {\tt HSL\_MA97} in parallel, for large problems in particular the cost
of the serial scaling can account for a much larger proportion of the runtime, making the reuse of
scaling an attractive option. It also suggests that, to avoid scaling being a bottleneck, parallel
scaling implementations are needed; this is something we plan to explore in the future.

Based on our findings, we have developed an Ipopt interface for the new
HSL linear solver {\tt HSL\_MA97} that offers the full range of scalings
tested in this paper. By including all the heuristics proposed in Section~\ref{dynamic scaling},
the interface provides users with a straightforward way to experiment to
find what works 
well for their problems. As we have observed that different scalings and
different heuristics work best for different problems, we would highly
recommend trying out the different options, particularly if many problems of a similar type are
 to be solved.

Finally, we note that the linear solver {\tt HSL\_MA97} together with the
scaling routines {\tt MC30}, {\tt MC64} and {\tt MC77}
and ordering/scaling routine {\tt HSL\_MC80}
are included in the 2011 release of HSL. The scaling routine
{\tt MC19} is part of the HSL Archive.
All use of HSL routines requires a licence; details of
how to obtain a licence and the routines are available at
{\tt http://www.hsl.rl.ac.uk/ipopt}. Note that routines in the HSL Archive are 
offered for free personal
use to all while those in HSL 2011 are available
without charge to academics for their teaching and research. In all cases, a separate commercial agreement is required to
allow redistribution of the code and/or the resulting binaries.

\vspace{0.5cm}
\noindent
{\large \bf Acknowledgements} \\ \\
We are grateful to our colleague Nick Gould for supplying the CUTEr test set and for commenting on a draft of this paper.

\bibliography{ref}

\end{document}